\documentclass[final,5p,times,twocolumn,authoryear]{elsarticle}

\usepackage{amssymb}
\usepackage{amsthm}
\usepackage{amsmath}
\usepackage{cleveref}
\usepackage{placeins}
\usepackage{float}

\newtheorem{rem}{Remark}

\journal{Computers and Chemical Engineering}

\begin{document}

\begin{frontmatter}



\title{Dynamic modeling and simulation of an electric flash clay calcination plant \\for green cement production}


\author[DTU]{Nicola Cantisani}

\author[DTU]{Jan Lorenz Svensen}
\author[FLS]{Shanmugam Perumal}
\author[DTU,CERE]{John Bagterp Jørgensen\corref{cor1}}
\ead{jbjo@dtu.dk}
\cortext[cor1]{Corresponding author}

\affiliation[DTU]{organization={Department of Applied Mathematics and Computer Science, Technical University of Denmark},
    city={Kongens Lyngby},
    postcode={2800}, 
    country={Denmark}}
    
\affiliation[CERE]{organization={Center for Energy Resources Engineering (CERE), Technical University of Denmark},
    city={Kongens Lyngby},
    postcode={2800}, 
    country={Denmark}}
    
\affiliation[FLS]{organization={FLSmidth Cement},
                 city={Chennai}, 
                 country={India}}
                 
\begin{abstract}
We present a novel dynamic model of an electric flash clay calcination plant. Calcined kaolinite-rich clay has been identified as one of the most effective candidates for supplementary cementitious material (SCM), because of its large availability. Calcination of clay is achieved via the dehydroxylation reaction, which does not release CO$_2$ (unlike limestone), but has a considerable energy requirement. The required high temperature can be met by electric resistive heating of the working gas in the plant, that can be powered by renewable energy. Therefore, CO$_2$-free calcination of clay can be achieved. Up to 50\% of the limestone-based clinker can be substituted by calcined clay (CC), making the cement more sustainable. We consider a plant that consists of gas-material cyclones that pre-heat the clay, a calciner, and a gas-recirculation system with electric heating of the gas. The model is formulated as a system of differential-algebraic equations (DAE). The model consists of thermophysical properties, reaction kinetics and stoichiometry, transport, mass and energy balances, and algebraic constraints. The model can be used to perform dynamic simulations with changing inputs, process design, and optimization. Moreover, it can be used to develop model-based control, which is relevant for flexible operation of a clay calcination plant for green cement production.
\end{abstract}




\end{frontmatter}



\section{Introduction}
Cement manufacturing is one of the largest sources of carbon dioxide emissions. It accounts for 8-9\% of the global carbon dioxide emission. Moreover, its manufacturing is responsible for 2-3\% of the global energy use, and its annual production is expected to grow by 50\% by 2050 \citep{Monteiro:etal:2017}. Accordingly, there is high interest in developing improved production methods that can reduce emissions. Clinker is the main component of cement, and it is produced by calcining limestone in the pyro-section (cyclone pre-heater tower, calciner, rotary kiln, and cooler) of a cement plant \citep{svensen2024firstengineering}. The calcination process of limestone is achieved via the endothermic reaction
\begin{equation}\label{eq:calcinationlimestone}
    \mathrm{CaCO_3} \rightarrow \mathrm{CaO + CO_2}.
\end{equation}
It is estimated that 1 kg of cement produces 0.6-0.8 kg of CO$_2$ \citep{Habert:2014}. Despite the number not being as high as other intensive industries (e.g. aluminium), the real environmental impact of the cement industry is due to its very high production volume. The annual global cement production has reached 2.8 billion tonnes, and is expected to reach 4 billion tonnes per year in 2050. Around 40\% of the CO$_2$ emissions are due to the burning of fossil fuel (typically coal) in the kiln. 50\% of the CO$_2$ emissions are related to the chemical process of calcination of limestone \eqref{eq:calcinationlimestone}, and the remaining 10\% of the CO$_2$ emissions are indirect emissions. Emissions reduction can be achieved in three different ways: 1) by lowering the clinker-to-cement ratio, 2) by substituting fossil fuel with renewable energy, 3) by capturing and storing the CO$_2$ generated by the process and fossil fuel combustion. The focus in this paper is related to the first strategy which implies that a part of the clinker is substituted by some supplementary cementitious material (SSM). In recent years, calcined kaolinite-rich clay (CC) as a clinker substitute has gained a lot of momentum, because of its abundance in nature and its CO$_2$-free calcination process. Substitution of up to 50\% clinker with CC in cement blends is viable, achieving similar mechanical properties and even improving some aspects of durability \citep{Scrivener:etal:2018,Habert:2014}. Limestone calcined clay cements are referred to as LC$^3$. Calcination of clay is achieved by thermally releasing the water bounded in the clay molecules.
By electrifying the clay calcination process using renewable energy and using the CC in the blend, emissions reduction of up to 50\% on the final cement product can be achieved. Furthermore, the use of electricity instead of fuel in the clay calcination process enables better temperature control, and thus higher product quality.

Because of the intermittent nature of renewable energy sources, it is relevant to be able to dynamically simulate and predict the effect of varying power input, as well as other process conditions (e.g. different clay compositions). A dynamic model of the process can serve this purpose. After fitting the parameters to the data recorded from a real plant, the dynamic model can become a Digital Twin of the plant, providing an inexpensive way to perform simulations without having to waste resources and time. Moreover, such a model can be exploited to perform process design and optimization. The effect of key parameters can be assessed and steady-state optimization can be performed. This can be used to look for the most cost-effective steady state to run the plant at. The nature (i.e. stable or unstable) of different steady states can also be discovered. Finally, the model can unlock the development of model-based control techniques, like model predictive control (MPC), for flexible and optimized process operation. Alternatively, the model can also be used to tune classical controllers (like PIDs). and compare performances of different controllers. It is therefore clear that such a tool has a high potential value, and can accelerate the development of the technology. In this paper, we develop such a simulation model for an electric flash clay calcination plant.

\cite{Laurini:etal:2024} present and discuss the vision on the electrification of clay calcination plants, and their integration into sustainable power grids. They present an Energy Management System (EMS) for a cement plant, that optimizes the power flows in the plant and handles the clay calcination plant as a smart load, providing flexibility. The dynamic model can then be used to simulate the effect and the plant output when a power set-point has been sent from the EMS.

There is very little literature available about dynamic modeling of the clay calcination process. \cite{Eskelinen:etal:2015} present a dynamic model for clay calcination in a multiple hearth furnace. They assume constant heat capacity for the solid phase and the model cannot be solved with a standard solver, but requires a special solution algorithm. \cite{Teklay:etal:2014} instead propose a transient one-dimensional model of the clay particles subject to calcination.

In this paper, we present a novel complete dynamic model of an industrial flash clay calcination plant. The plant consists of gas-solid cyclones that pre-heat the clay, a calciner, and a gas-recirculation system with electric heating of the gas. The model is formulated as a system of differential-algebraic equations (DAE). The model is based on first engineering principles, i.e. mass and energy balances. We use a rigorous approach that incorporates thermodynamic functions as algebraic constraints. This allows us to handle complicated (non-constant) expressions of the heat capacity, and have the state variables, such as temperature and pressure, as algebraic variables. This technique makes the model not only more realistic, but easy to implement as a standard DAE system. Moreover, our formulation in blocks allows easy modifications, if needed. The main physical units that we model in detail are the calciner and the gas-solid cyclones.

This paper is an extension of a conference paper \citep{Cantisani:etal:2024}. \cite{Cantisani:etal:2024} model the calciner and provide the chemical and thermophysical models. In this paper, we model the entire entire clay calcination plant using the calciner model as well as the chemical and thermophysical model provided by \cite{Cantisani:etal:2024}.

\subsection{Paper outline}
The paper is structured as follows. Section \ref{sec:process_description} provides a short overview and a description of the clay calcination plant. Section \ref{sec:model} presents the mathematical model of the process, with all of its units and the connections. Each subsection describes mathematically each unit, and finally an overview of the model (i.e. its equations and variables) is given. The Jacobian of the model is also discussed. Section \ref{sec:simulation} presents some simulation results of the model. Section \ref{sec:conclusion} concludes the paper.

\section{Process description}\label{sec:process_description}
\begin{figure}[tb]
    \centering
    \includegraphics[width=0.5\textwidth]{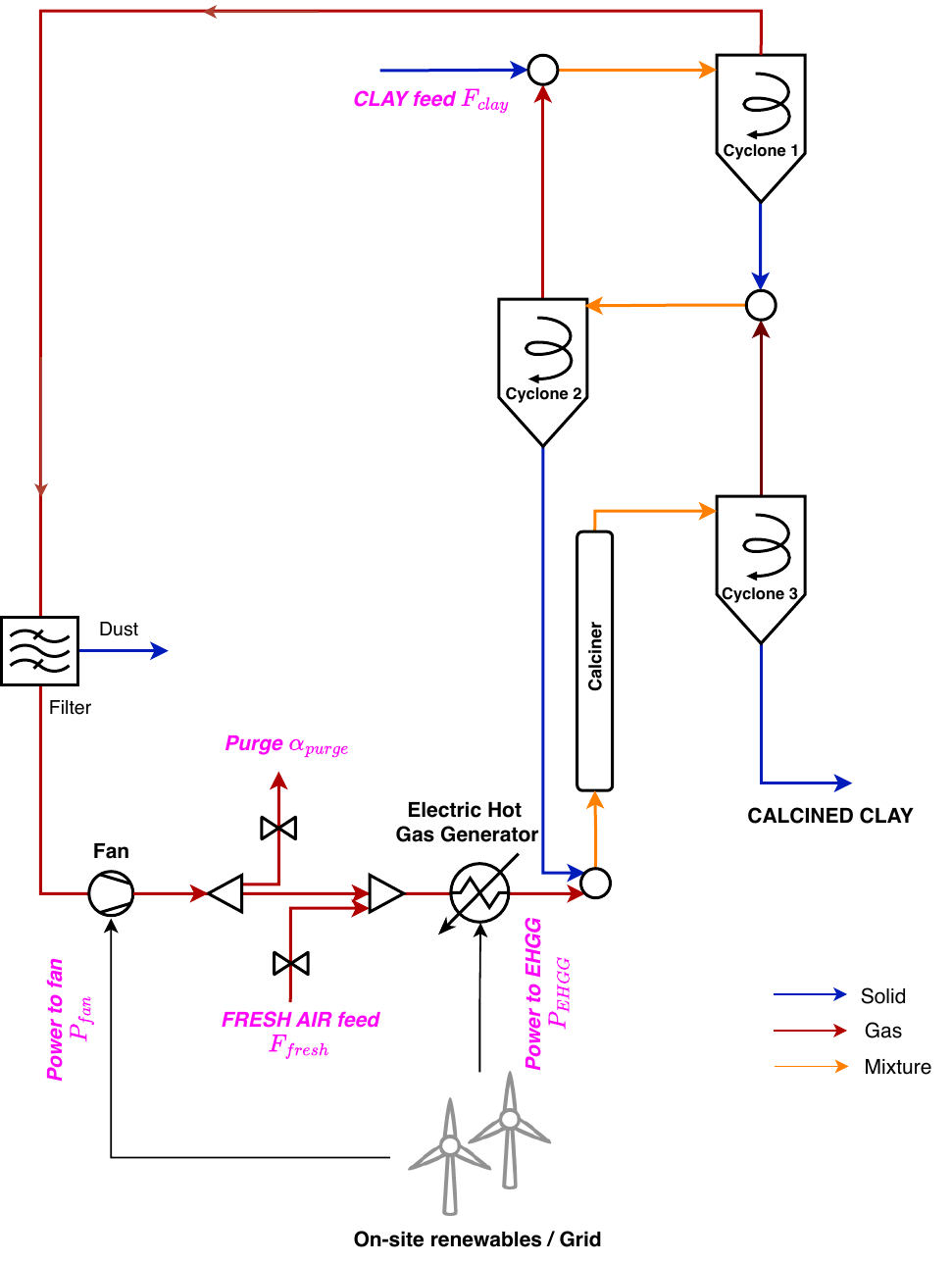}
    \caption{Process diagram of the electric flash clay calcination. The fresh clay gets pre-heated through two cyclones, before reaching the calciner, where the calcination reaction is completed. The circulating gas is heated with renewable power and recycled in the loop. The purge ensures that water does not accumulates in the loop. The main purpose of the process is to remove the water bound in the clay to get calcined clay (water free clay). The system inputs are highlighted in magenta.}
    \label{fig:processdiagram}
\end{figure}
Figure \ref{fig:processdiagram} shows a diagram of the clay calcination process that we consider. In this section, we provide a short description of the clay calcination process.
The thermal activation of the clay is performed in a loop. The fresh clay is introduced in the loop after being crushed, at the inlet of Cyclone 1 (see Figure \ref{fig:processdiagram}). The material undergoes pre-heating through two cyclones (Cyclone 1 and Cyclone 2), where a part of the clay already gets calcined because of the high temperature. The pre-heated solid is then introduced in the calciner. The calciner is a long plug-flow reactor (PFR) where the solid material stream is mixed with the hot gas stream coming from the electric hot gas generator. The hot gas transfers heat to the solid particles, ensuring that all the clay gets calcined. The last cyclone (cyclone 3 in Figure \ref{fig:processdiagram}) separates the solid product from the gas before leaving the process. The gas is recirculated in a loop, in order to recover energy. A circulating fan holds a certain pressure difference in the whole loop. A ceramic filter removes any solid particle residue (dust) from the gas before being recycled. Some of the gas is purged to remove water and ensure that water does not accumulates in the loop. Fresh air is mixed to the recycled stream. The gas flow undergoes heating by passing through the electric hot gas generator, which transfers heat to the gas by resistive heating. This could, alternatively, be replaced by a thermal storage, that works as a buffer when electricity is intermittent. The process is CO$_2$ free, as the hot gas generator runs on renewable energy and the clay does not release any carbon dioxide.

\begin{figure}
    \centering
    \includegraphics[width=0.48\textwidth]{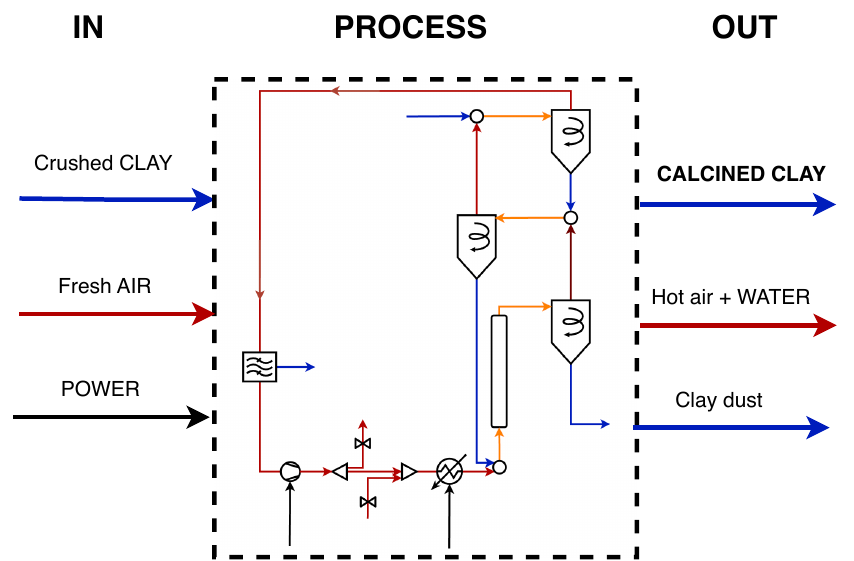}
    \caption{Overview of the inputs and outputs of the process.}
    \label{fig:InputOutput}
\end{figure}

Figure \ref{fig:InputOutput} provides a visual overview of the inputs and outputs of the process. The inputs are the clay (flow rate, composition, temperature), the fresh air, and the power. The outputs are the calcined clay (flow rate, composition, temperature), the hot air and water in the purge stream, and the clay dust. 
\section{Mathematical model description}\label{sec:model}
The modeling approach adopted in this paper is based on establishing first-principles equations for the single units of interest in the process. The units are therefore modeled \textit{dynamically}, and they are connected appropriately. Additional secondary units, that are not of particular interested, may be modeled by steady state relationships. 
By collecting all the equations, a system of differential-algebraic equations arises in the form
\begin{subequations}
\label{eq:index1DAEsystem}
    \begin{align}
        \frac{d x}{dt} & = f(x,y,u,d,p), \\
        0 & = g(x,y,u,d,p).
    \end{align}
\end{subequations}
$x$ are the differential variables, $y$ are the algebraic variables, $u$ are the inputs (or manipulated variables), $d$ are the disturbances, and $p$ the model parameters. The time dependency is not included explicitly, as the model is time invariant.
\subsection{Model structure}
The model is an index-1 differential algebraic equation system \eqref{eq:index1DAEsystem} that consists of several building blocks. These elements are obtained using a systematic combination of stoichiometry and kinetics, thermophysical model functions, transport functions, mass- and energy-balances, and equipment related equations that provide the algebraic equations. 

The building blocks for the systematic model constructions can be summarized as
\begin{itemize}
    \item Stoichiometry and kinetics of the dehydroxilation reaction of the clay.
    \item Thermophysical model (enthalpy, volume and internal energy of the materials, and mixture viscosity). 
    \item Model of the main units (dynamic),
    \item Model of the secondary components (static),
    \item Connection of the units.
\end{itemize}
The chemical and thermophysical models are shared by all the units. We model extensively (and therefore dynamically) the calciner and the cyclones. Each of these two models consists of the following subparts
\begin{itemize}
    \item Mass balance,
    \item Energy balance,
    \item Transport model (pressure drop and velocities, material fluxes),
    \item Algebraic relations related to the thermophysical properties and dimension of the unit.
\end{itemize}
The same modeling technique is used by \cite{Rosbo:etal:2023} and \cite{martinsen2023a} for an ammonia synthesis loop.

\subsection{Notation}
We specify all the concentrations, $c$, with respect to the reactor volume, $V$, i.e.
\begin{equation}\label{eq:conc_def}
    c = \frac{n}{V}.
\end{equation}
$n$ indicates number of moles. The following notation for the materials is used throughout the paper
\begin{itemize}
    \item AB$_2$ = kaolinite (solid),
    \item A = metakaolin (solid),
    \item B = water (gas),
    \item air = dry air (gas),
    \item Q = quartz (solid). 
\end{itemize}
The concentrations, $c$, the molar fluxes, $N$, the number of moles, $n$, and the molar masses, $M$, are always intended as vectors with the following order
\begin{equation}
    c = \begin{bmatrix}
        c_{AB_2}\\ c_{A} \\ c_{B} \\ c_{air} \\ c_{Q}
    \end{bmatrix}, \quad N = \begin{bmatrix}
        N_{AB_2}\\ N_{A} \\ N_{B} \\ N_{air} \\ N_{Q}
    \end{bmatrix}, \quad n = \begin{bmatrix}
        n_{AB_2}\\ n_{A} \\ n_{B} \\ n_{air} \\ n_{Q}
    \end{bmatrix}, \quad M = \begin{bmatrix}
        M_{AB_2}\\ M_{A} \\ M_{B} \\ M_{air} \\ M_{Q}
    \end{bmatrix}.
\end{equation}
The subscript $s$ and $g$ indicate the vector of only solid phase and gas phase, e.g.
\begin{equation}
c_s = \begin{bmatrix}
        c_{AB_2}\\ c_{A} \\ c_{Q}
    \end{bmatrix}, \quad c_g = \begin{bmatrix}
        c_{B} \\ c_{air}
    \end{bmatrix}.
\end{equation}
The reader should be aware that, when possible, the naming of the variables is generally related to its section only, unless otherwise specified. That means that, for example, the concentrations variable $c$ in the calciner section is only referred to the calciner. This is clear from the context, and it is done in order to keep the notation as simple as possible. Should there be any ambiguity, or variables from different units are combined, the superscripts ${\cdot}^{Calc}$,${\cdot}^{Cyc1}$,${\cdot}^{Cyc2}$,${\cdot}^{Cyc3}$ are used. These indicate the calciner, cyclone 1, cyclone 2 and cyclone 3, respectively, as defined in Figure \ref{fig:processdiagram}.
\subsection{Stoichiometry and kinetics}
We consider clay that is composed of kaolinite and quartz. The main reaction occurring when calcinating clay is dehydroxylation of kaolinite. The reaction leads to the formation of metakaolin and water vapor, according to the following reaction
\begin{equation}
        \mathrm{Al_2O_2 \cdot 2SiO_2 \cdot 2H_2O (s) \rightarrow } \mathrm{ Al_2O_2 \cdot 2SiO_2 (s)} \mathrm{ + 2 H_2O (g).}
\end{equation}
The temperature range of the reaction is 450-700 ${}^\circ$C. The reaction is endothermic. Notice that quartz is not involved in the chemical reaction (it is therefore inert). Using our notation, the reaction is expressed as
\begin{equation}
\mathrm{AB_2 \rightarrow A + 2 B}.
\end{equation}
We model the reaction kinetics as a third-order reaction, with activation energy $E_A=202$ kJ/mol and pre-exponential factor $k_0 = 2.9 \times 10^{15}$ s$^{-1}$ \citep{Ptacek:etal:2010}. The reaction rate may therefore be expressed as
\begin{equation} \label{eq:reactionrate}
    r = r(c_{AB_2},T_s) = k \, c_{AB_2}^3
\end{equation}
where 
\begin{equation}
    k = k(T_s) = k_0 \exp \left( - \frac{E_A}{R_{gas} T_s} \right).
\end{equation}
$T_s$ is the temperature of the solid. The chemical production rate is
\begin{equation}
    R = \nu'r(c),
\end{equation}
where $\nu$ is the stoichiometric matrix.
The stochiomatric matrix has the form
\begin{equation}
     \nu = [-1,1,2,0,0].
\end{equation}
Notice that air and quartz do not participate in the reaction, but they are involved in the mass transport and heat exchange.

\subsection{Thermophysical model}
Expressing energy balances in terms of internal energy makes the model depend on temperature $T$, pressure $P$, and number of moles $n$. By using a thermodynamic library, or by constructing one, we can compute the volume $V$, the enthalpy $H$, and (therefore) the internal energy $U$, i.e.
\begin{subequations}
    \begin{align}
        V & = V(T,P,n),  \\
        H &= H(T,P,n), \\
        U &= H - PV.
    \end{align}
\end{subequations}
Moreover, notice that the volume and enthalpy functions of a mixture may be computed as
\begin{subequations} \label{eq:mixture}
    \begin{align}
        V(T,P,n) & = \sum_{i} n_i v_i(T,P) = n^\top v(T,P),  \\
        H(T,P,n) & = \sum_{i} n_i h_i(T,P)= n^\top h(T,P),
    \end{align}
\end{subequations}
where the index $i$ indicate the $i$-th species of the mixture. $h$ and $v$ indicate molar enthalpy and volume. \eqref{eq:mixture} is only valid when the molar enthalpies do not depend on the molar fractions, like in our case. Since these functions are homogeneous of order 1 with respect to the number of moles, we can divide by the reactor volume and obtain the volumetric quantities. 
We use the $\hat \cdot$ notation to indicate them.
\begin{subequations}
    \begin{align}
        \hat v & = V(T,P,c) = 1, \\
        \hat h &= H(T,P,c), \\
        \hat u &= \hat h - P \hat v.
    \end{align}
\end{subequations}
In the same way, we may compute the enthalpy flux (indicated by $\Tilde \cdot$)
\begin{subequations}
    \begin{align}
    \Tilde H = H(T, P, N).
    \end{align}
\end{subequations}


\subsubsection{Solid phase}
The solid phase is kaolinite, metakaolin, and quartz.
The molar enthalpy of each component $i$ may be evaluated if the expression of the heat capacity $c_{p,i}$ of the material is available.
\begin{subequations}\label{eq:molar_enthalpy}
    \begin{align}
        h_i(T,P) & = h_i(T^0,P^0) + \int_{T^0}^T c_{p,i}(s) ds
    \end{align}
\end{subequations}
We use the enthalpy of formation of the material at standard conditions for the reference state. That is
\begin{equation}
    h_i(T^0=298.15 \text{ K},P^0=1 \text{ bar}) = \Delta H^0_{form,i}.
\end{equation}
Table \ref{tab:heatcap} reports the values for metakaolin and kaolinite. We use the notation
\begin{subequations}
\begin{align}
        H_s & = H_s(T_s,P,n_s), \\
        V_s & = V_s(T_s,P,n_s),
\end{align}
\end{subequations}
to indicate the enthalpy and the volume of the solid material.

The heat capacity of the solid components may given by the following empirical expression 
\begin{equation}
    c_{P}(T) = k_1 + k_2 T + k_3 T^2 + \frac{k_4}{T} + \frac{k_5}{T^2} + \frac{k_6}{\sqrt{T}},
\end{equation}
with $\{k_i\}_{i=1,...,6}$ being coefficients.
This is only valid within the temperature interval $T_{min} \leq T \leq T_{max}$. Table \ref{tab:heatcap} reports the coefficients for metakaolin and kaolinite \citep{Gariz:2023,Bale:etal:2016}. The temperature $T$ must be given in K and the heat capacity $c_P$ is in 
$\mathrm{\frac{J}{mol \cdot K}}$.
The molar volume is given by
\begin{equation}
    v(T) = v_1 + v_2 T.
\end{equation}
Table \ref{tab:volume} reports $v_1$ and $v_2$ for metakaolin and kaolinite \citep{Gariz:2023,Bale:etal:2016}. \cite{Lindstrom:2023} provide the enthalpy and the volume expressions for quartz.

\begin{table}[tb]
    \centering
    \caption{Coefficients for the solid heat capacity and standard enthalpy of formation of metakaolin and kaolinite.}
    \begin{tabular}{|c|c|c|}
    \hline
        \textbf{Coefficients} & \textbf{Metakaolin} & \textbf{Kaolinite}  \\
        \hline
        $k_1$ & 2.294924 $\times 10^{2}$ $\mathrm{\frac{J}{mol \cdot K}}$ & 1.4303 $\times 10^{3}$ $\mathrm{\frac{J}{mol \cdot K}}$ \\
        $k_2$ & 3.68192 $\times 10^{-2}$ & -7.886 $\times 10^{-1}$ \\
        $k_3$ & 0 & 3.034 $\times 10^{-4}$\\
        $k_4$ & 0 & 0 \\
        $k_5$ & -1.456032 $\times 10^{6}$ & 8.334 $\times 10^{6}$ \\
        $k_6$ & 0 & -1.862 $\times 10^{4}$ \\
        \hline
        $T_{min}$ & 298 K & 298 K\\
        $T_{max}$ & 1800 K & 700 K \\
        \hline
        $\Delta H^0_{form}$ & -3.211 $\times 10^{6}$ $\mathrm{\frac{J}{mol}}$ & -4.11959 $\times 10^{6}$ $\mathrm{\frac{J}{mol}}$ \\
        \hline
    \end{tabular}
    \label{tab:heatcap}
\end{table}
    
\begin{table}[tb]
    \centering
    \caption{Coefficients for the molar volume.}
    \begin{tabular}{|c|c|c|c|}
    \hline
        \textbf{Coefficients} & \textbf{Metakaolin} & \textbf{Kaolinite} \\
        \hline
        $v_1$ & 41.4736 $\mathrm{\frac{m^3}{mol}}$ & 30 $\mathrm{\frac{m^3}{mol}}$\\
        $v_2$ & 3.39116 $\times 10^{-3} \mathrm{\frac{m^3}{mol \cdot K}}$  & 0 $\mathrm{\frac{m^3}{mol \cdot K}}$\\
        \hline
    \end{tabular}
    \label{tab:volume}
\end{table}

\subsubsection{Gas phase}
The gas phase in the clay calcination process consists of dry air and water vapor. We treat air as an ideal gas mixture made of 78\% nitrogen (N$_2$), 21\% oxygen (O$_2$), and 1\% argon (Ar). The molar enthalpy of each component is computed as for the solid \eqref{eq:molar_enthalpy}. We indicate the enthalpy and the volume of the gas mixture with 
\begin{subequations}
\begin{align}
        H_g & = H_g(T_g,P,n_g), \\
        V_g &= V_g(T_g,P,n_g).
\end{align}
\end{subequations}
$T_g$ indicates the temperature of the gas phase. 
\cite{Lindstrom:2023} provide the enthalpy expressions for nitrogen, oxygen, argon, and water vapor.

The molar volume of a single gas may be computed using the ideal gas law
\begin{equation}
    v_{g,i} =  \frac{R T_g}{P}.
\end{equation}

\subsubsection{Viscosity}
The viscosity of a suspended gas mixture may be computed by the extended Einstein equation of viscosity \citep{Kiyoshi:2006,svensen2024firstengineering}.
\begin{equation}
    \mu = \mu_g \frac{1 + \hat v_s/2}{1 - 2 \hat v_s}.
\end{equation}
$\mu_g$ is the viscosity of the gas phase. The viscosity of a mixed gas may be given as 
\begin{subequations}
\begin{align}
\mu_g & =\sum_i \frac{x_i \, \mu_{g, i}}{\sum_j x_j \phi_{i j}}, \\
\phi_{i j} & =\left(1+\sqrt{\frac{\mu_{g, i}}{\mu_{g, j}}} \sqrt[4]{\frac{M_j}{M_i}}\right)^2\left(2 \sqrt{2} \sqrt{1+\frac{M_i}{M_j}}\right)^{-1}.
\end{align}
\end{subequations}
$\mu_{g,i}$ is the viscosity of a single gas \citep{Wilke:2004}. $x_i$ and $M_i$ indicate the molar fraction and the molar mass of the $i$-th gas, respectively. In our case $i = \{B,air\}$. The viscosity of a single gas can be expressed as a function of the temperature \citep{Sutherland:1893}
\begin{equation}
    \mu_{g,i} = \mu_0 \Big(\frac{T_g}{T_0} \Big)^\frac{3}{2} \frac{T_0 + S_\mu}{T_g + S_\mu}.
\end{equation}
$S_\mu$ can be calibrated given two measures of viscosity.

\section{Calciner}
This subsection presents the partial differential equations (PDAE) model of the calciner \citep{Cantisani:etal:2024}. The calciner is modeled as a PFR of length $L$ and diameter $d$. The model consists of mass and energy balances in space and time, transport equations for the velocities and the fluxes, and algebraic relations. 
The model is discretized in space and converted to DAEs. 

\subsection{Transport model}
Mass and energy balances in the calciner directly depend on the spatial material flux. We consider advection and diffusion for the flux. The flux depends directly on the flow velocity.

\subsubsection{Velocity} 
The velocity along the reactor, $v$, may be modeled as a function of the pressure drop, $\Delta P$, along the length, $\Delta z$, using the Darcy-Weisbach equation for turbulent flows \citep{svensen2024firstengineering}.
\begin{equation}
    v = v\Big( \frac{\Delta P}{\Delta z} \Big) = \Big( \frac{2}{0.316} \sqrt[4]{\frac{d^5}{\mu \, \rho^3}} \frac{|\Delta P|}{\Delta z} \Big)^{ \frac{4}{7} } \mathrm{sgn} \Big( \frac{\Delta P}{\Delta z} \Big).
\end{equation}
$\rho$ is the density of the mixture. The equation is only valid for Mach number $<$ 0.2. The density of the mixture may be computed as
\begin{equation}
    \rho = M \cdot c 
\end{equation}

\subsubsection{Advection and diffusion}
The molar flux of the mixture, $N$, is modeled as the sum of advection, $N_a$, and Fick's diffusion, $N_d$.
\begin{equation}
    N = N_a + N_d,
\end{equation}
where
\begin{subequations}
    \begin{align}
         N_a & = v \cdot c, \\
         N_d & = -D  \odot \partial_z c.
    \end{align}
\end{subequations}
$D$ are the diffusion coefficients, and $\odot$ is the element-wise product. 

\subsection{Mass balance} 
The mass balance in the calciner reads (PFR model in \citep{Nielsen:2023})
\begin{equation}\label{eq:massbalance}
    \partial_t c = - \partial_z N + R.
\end{equation}
The PDE \eqref{eq:massbalance} describes the dynamics of the concentrations $c(t,z)$ in time and space. The variable $z$ is the calciner length direction. $N(t,z)$ are the fluxes and $R(c)$ is the chemical production rate.

\subsection{Energy balance}
Because two different phases (solid and gas) are interacting and exchanging energy in the calciner, we need to keep track of the temperature of both of them. We hereby derive an energy balance for the solid and the gas phases.

Let us consider a control volume $\Delta V = A \Delta z$. Let us assume that the solid particles are well mixed with the gas, and that they have identical shape and size. We assume that each solid particle is a perfect ball of radius $r_b$.
The accumulated energy \textit{in the solid phase} in the control volume, during the time interval $\Delta t$, is
\begin{equation}\label{eq:energy1}
\begin{split}
    \Delta U_s =& A \, N_s(t,z) h_s(T_s,P) \Delta t - A N_s(t,z+\Delta z) h_s(T_s,P) \Delta t  \\ & + J_{sg} \Delta t - Q_{amb}\Delta t.
\end{split}
\end{equation}
$J_{sg}$ is the heat transfer rate between solid and gas phases, $Q_{amb}$ is the heat loss to the ambient. We can express the heat transfer between the phases as
\begin{equation} \label{eq:sg}
   J_{sg} = k_{sg} A_{sg} (T_g - T_s).
\end{equation}
$k_{sg}$ and $A_{sg}$ are the solid-to-gas heat transfer coefficient and the transfer area, respectively.
Because of the assumption on the solid particles, we can compute the transfer area between the solid and the gas as
\begin{equation}
    A_{sg} = n_{ball} A_{ball} = \frac{V_s}{V_{ball}}A_{ball} = \frac{3 V_s}{r_b}.
\end{equation}
$n_{ball}$ is the number of balls in the volume, $V_{ball}$ is the volume of a ball and $A_{ball}$ is the area of a ball.

Now inserting \eqref{eq:sg} in \eqref{eq:energy1}, and dividing by $\Delta V$ and $\Delta t$, we get 
\begin{equation}
    \begin{split}
         \frac{\Delta U_s}{\Delta t \Delta V} = & - \frac{N_s(t,z+\Delta z) h_s(T_s,P) - N_s(t,z) h_s(T_s,P)}{\Delta z} \\
    & + k_{sg} \frac{3 V_s}{r_b \Delta V} (T_g - T_s) - \frac{Q_{amb}}{\Delta V} 
    \end{split}
\end{equation}

Using volumetric quantities and letting $\Delta t \rightarrow 0$ and $\Delta z \rightarrow 0$, the following PDE arises
\begin{equation}
    \frac{\partial \hat u_s}{\partial t} = - \frac{\partial \tilde H_s}{\partial z}  + k_{sg} \frac{3 \hat v_s}{r_b} (T_g - T_s) - \hat Q_{amb}.
\end{equation}

The same derivation can be repeated for the gas phase. In compact notation, the following set of PDE describes the (volumetric) energy balance of the solid and gas phases in the calciner
\begin{subequations}\label{eq:energybalance}
    \begin{align}
        \partial_t \hat u_s & = - \partial_z \Tilde H_s  + \hat J_{sg} - \hat Q_{amb}, \\
        \partial_t \hat u_g & = - \partial_z \Tilde H_g - \hat J_{sg} - \hat Q_{amb}.
    \end{align}
\end{subequations}
In reality, the solid particles have different sizes. If the particle size distribution is available, then the median diameter $d_{med}$ is used in the model. Therefore we rather express the the volumetric heat transfer between solid and gas as
\begin{equation}\label{eq:solidtogas}
    \hat J_{sg} = k_{sg} \frac{6 \hat v_s}{d_{med}} (T_g - T_s).
\end{equation}

\subsection{Algebraic relations}
Some extra algebraic equations are needed to solve the system. The volume of solid and gas phases should sum to the reactor volume, i.e.
\begin{equation} \label{eq:volume_alg}
    V_s(T_s,P,c_s) + V_g(T_g,P,c_g) - 1 = 0.
\end{equation}
Moreover, the differential variables $\hat u_s$ and $\hat u_g$ should match the quantities computed via the thermodynamic functions, that is
\begin{subequations} \label{eq:intenergy_alg}
    \begin{align}
         & U_s(T_{s},P,c_{s}) - \hat u_{s} =0, \\
         & U_g(T_{g},P,c_{g}) - \hat u_{g} =0.
    \end{align}
\end{subequations}

\subsection{PDE spatial discretization} \label{sec:discretization}
\begin{figure}[tb]
    \centering
    \includegraphics[width=0.48\textwidth]{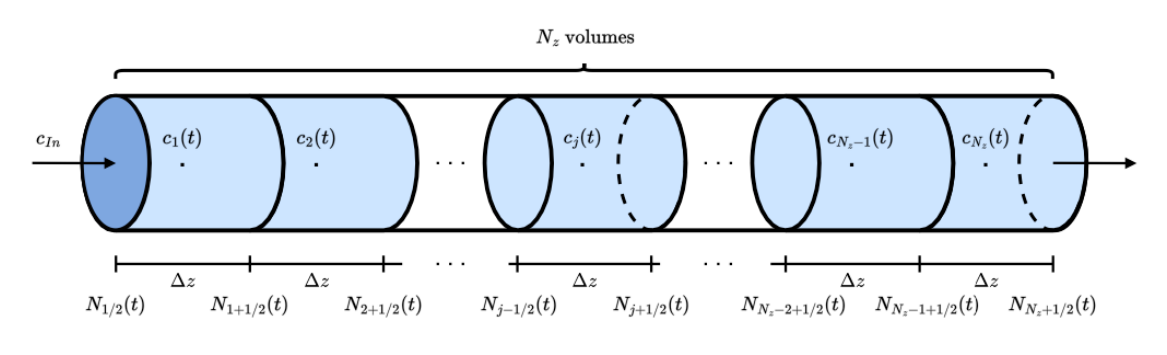}
    \caption{Finite volume discretization of the calciner.}
    \label{fig:finitevolume}
\end{figure}
The mass and energy balances for the calciner resulted in a set of partial differential equations in space and time, \eqref{eq:massbalance} and \eqref{eq:energybalance}. In order to be solved and simulated, we perform spatial discretization by dividing the calciner into finite volumes along the length $z$. We apply central difference approximation to evaluate the derivatives at the center of each cell. Let us consider $N_z$ finite volumes (see Figure \ref{fig:finitevolume}). The fluxes at the cell interfaces are
\begin{equation}
    N_{i+1 / 2} =v_{i+1/2} \, c_i-D \odot \frac{c_{i+1}-c_i}{\Delta z},
\end{equation}
where the velocities are such that
\begin{equation}
    v_{i+1/2} = v\Big(\frac{P_{i+1}-P_{i}}{\Delta z}\Big),
\end{equation}
for $i \in\left\{1,2, \ldots, N_z-1\right\}$.
At the first and the last cell interface we have
\begin{subequations}
    \begin{align}
        N_{1/2} & = v_{1/2} c_{in},\\
        N_{N_z+1/2} & = v_{N_z+1/2} c_{N_z}.
    \end{align}
\end{subequations}
$c_{in}$ are the inlet concentrations.
The velocities at the inlet and outlet of the calciner, namely $v_{1/2}$ and $v_{N_z+1/2}$, depend on the pressures before the inlet and after the outlet, $P_{in}$ and $P_{out}$, respectively.
\begin{subequations}
\begin{align}
    v_{1/2} &= v\Big(\frac{P_1-P_{in}}{\Delta z}\Big), \\
    v_{N_z+1/2} &= v\Big(\frac{P_{out}-P_{N_z}}{\Delta z} \Big).
\end{align}
\end{subequations}

The spatial discretization of the mass balance equation results in $N_z$ ordinary differential equations (ODE)
\begin{equation} \label{eq:massODE}
    \frac{d c_i}{d t}=-\frac{N_{i+1 / 2}-N_{i-1 / 2}}{\Delta z}+R\left(c_i\right),
\end{equation}
for $i \in\left\{1,2, \ldots, N_z\right\}$. Notice that each $c_i$ is a vector with the 5 chemical components. Therefore the mass balance results in $5 \cdot N_z$ equations.

In the same way, the energy balances are also discretized in space
\begin{subequations} \label{eq:energyODE}
\begin{align}
     \frac{d \hat u_{s,i}}{dt} & = - \frac{\tilde H_{s,i+1/2}-\tilde H_{s,i-1/2}}{\Delta z} + \hat J_{sg,i} - \hat Q_{amb}, \\
 \frac{d \hat u_{g,i}}{dt} & = - \frac{\tilde H_{g,i+1/2}-\tilde H_{g,i-1/2}}{\Delta z} - \hat J_{sg,i} - \hat Q_{amb}.
\end{align}
\end{subequations}
The enthalpy fluxes of the solid material at the cell interfaces are
\begin{subequations}
    \begin{align}
    \tilde H_{s,1/2} &= H(T_{s,in},P_{in},N_{s,1/2}), \\
\tilde H_{s,i+1/2} &= H(T_{s,i},P_i,N_{s,i+1/2}),
    \end{align}
\end{subequations}
for $i \in\left\{1,2, \ldots, N_z\right\}$.
The enthalpy fluxes of the gas at the cell interfaces are
\begin{subequations}
    \begin{align}
    \tilde H_{g,1/2} &= H(T_{g,in},P_{in},N_{g,1/2}), \\
\tilde H_{g,i+1/2} &= H(T_{g,i},P_i,N_{g,i+1/2}),
    \end{align}
\end{subequations}
for $i \in\left\{1,2, \ldots, N_z\right\}$. The inlet temperatures of the solid and gas, $T_{s,in}$ and $T_{g,in}$, are manipulated variables.

The heat transfer term between solid and gas is
\begin{equation}
    \hat J_{sg,i} = k_{sg} \frac{6 \hat v_{s,i}}{d_{med}} (T_{g,i} - T_{s,i}).
\end{equation}
for $i \in\left\{0,1, \ldots, N_z\right\}$.

Finally, the algebraic equations are 
\begin{subequations} \label{eq:algebraic}
    \begin{align}
        & U_s(T_{s,i},P_i,c_{s,i}) - \hat u_{s,i} =0,\\
        & U_g(T_{g,i},P_i,c_{g,i}) - \hat u_{g,i} =0,\\
        & V_s(T_{s,i},P_i,c_{s,i}) +  V_g(T_{g,i},P_i,c_{g,i}) - 1 = 0, 
    \end{align}
\end{subequations}
for $i \in\left\{0,1, \ldots, N_z\right\}$.

\subsection{Summary}
The total discretized DAE model consists of \eqref{eq:massODE}, \eqref{eq:energyODE}, and \eqref{eq:algebraic}. This results in $10 \cdot N_z$ equations. The differential variables $x$ and the algebraic variables $y$ are
\begin{equation}
    x^{Calc} = \begin{bmatrix} c_i \\ \hat u_{s,i} \\ \hat u_{g,i} \end{bmatrix}, \quad y^{Calc} = \begin{bmatrix} T_{s,i} \\ T_{g,i} \\ P_{i} \end{bmatrix} 
\end{equation}
for $i \in\left\{1,2, \ldots, N_z\right\}$.

\section{Cyclones}

\subsection{Working principle and geometry} 
\begin{figure}[tb]
    \centering
    \includegraphics[width=0.45\textwidth]{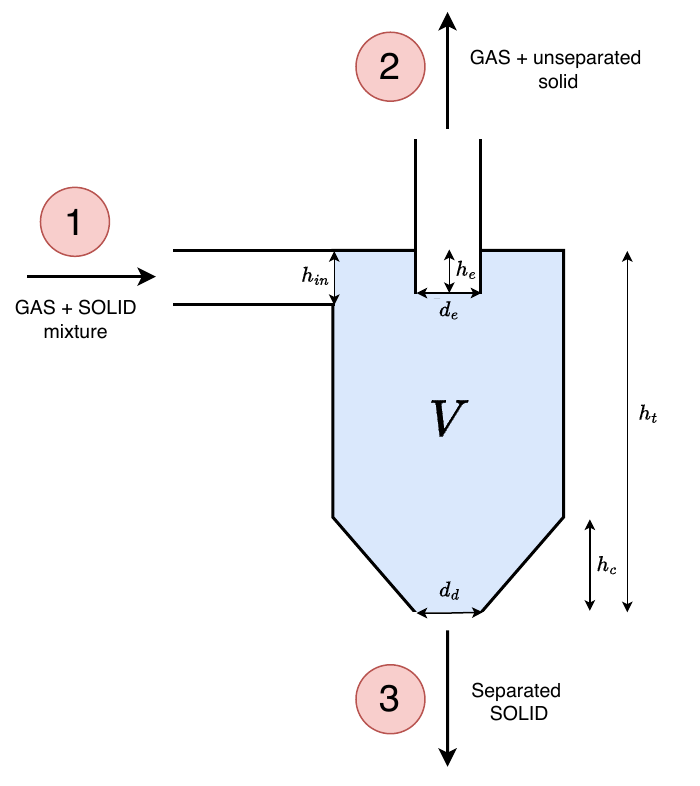}
    \caption{Cyclone section: geometric parameters definition, cell volume (in blue) and material streams numbering (in red).}
    \label{fig:cyclone_volume}
\end{figure}

\begin{figure}[tb]
    \centering
    \includegraphics[width=0.2\textwidth]{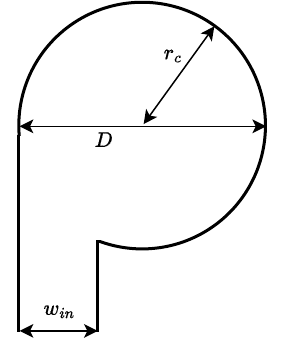}
    \caption{Axial profile of a cyclone with slot inlet, with geometric parameters definition.}
    \label{fig:cyclone_axial}
\end{figure}

A cyclone is a physical device that is used to separate suspended solid particles in a gas ('dedusting'). In the clay calcination process, they are also used to pre-heat the clay before it enters the calciner (cyclone 1 and 2 in Figure \ref{fig:processdiagram}). Figure \ref{fig:cyclone_volume} shows the section drawing of a cyclone. The numbering of the streams is assigned in the figure and used consistently as subscripts in this section, for fluxes, velocities, and cross sectional areas. Gas and solid enter the cyclone in the inlet (1). The mixture swirls in the cyclone body until it reaches the bottom where most of the solid drops in the separation outlet (3). The gas stream with the remaining unseparated solid particles is captured by the gas outlet (2) (or 'vortex finder'). Figure \ref{fig:cyclone_volume} also shows the geometric parameters definition, and the total cyclone volume, $V$, in blue. \cite{Hoffmann:2002} provide further description of the working principle.
Figure \ref{fig:cyclone_axial} shows the axial profile of the cyclone. We consider cyclones with slot (or 'tangential') inlet, i.e. rectangular section inlet, that ends tangentially to the cyclone body.

The total volume of the cyclone body can be computed as
\begin{equation}
    V = \pi  \left(r_c^2 (h_t - h_c) + \frac{h_c}{3} (r_c^2 + \frac{d_d^2}{4} + r_c \frac{d_d}{2}) - \frac{d_e^2}{4} h_e \right)
\end{equation}
The total surface area is
\begin{equation}
    A_c = 2 \pi r_c (h_t - h_c) + \pi \left( r_c^2 - \frac{d_e^2}{4} \right) + \pi \left( r_c + \frac{d_d}{2} \right) \sqrt{\left( r_c - \frac{d_d}{2} \right)^2 + h_c^2} 
\end{equation}
Beware that a collecting hopper under the cyclone is not included in our model, since the cyclones are connected in the process and the material is not collected until the end of the cycle.

\subsection{Model outline} 
As opposed to the calciner, we model the cyclones using a lumped approach. This means that we consider them as a single cell unit, with total volume $V$. Derivatives in space of the model variables do not appear. 

\subsection{Mass balance} 
The mass balance in the cyclone body is 
\begin{equation}\nonumber
    \frac{dn}{dt} = f_{1}-f_{2}-f_{3} + V R,
\end{equation}
where $f$ indicates a mass flow rate. By dividing by the total volume we obtain the concentration based mass balances
\begin{equation}\label{eq:cyclonemassbalance}
    \frac{dc}{dt} = \frac{1}{V} \Big( A_1 N_{1}-A_{2}N_{2}-A_{3}N_{3} \Big) + R,
\end{equation}
where $N$ indicates a flux and $A$ a cross-sectional area.

\subsection{Energy balance}
In the same way, the energy balances in the cyclone are derived. The accumulated internal energy in the solid phase over the time interval $\Delta t$ is
\begin{equation}\nonumber
\begin{split}
    \Delta U_s = &\Big( A_1 \, N_{1,s} h_s(T_{1,s},P_{in})- A_{2} \, N_{2,s} h_s(T_s,P) \\ & - A_{3} \, N_{3,s} h_s(T_s,P) \Big) \Delta t + J_{sg} \Delta t - Q_{amb}\Delta t.
\end{split}
\end{equation}
$T_{1,s}$ is the inlet solid temperature. Notice that inlet solid and gas have different temperatures.
By dividing by $V$ and letting $\Delta t \rightarrow 0$, and repeating the derivation for the gas phase, we obtain the energy balance for the solid and gas phases in the cyclone
\begin{subequations}
\label{eq:cycloneenergybalance}
\begin{align}
    \frac{d \hat u_s}{dt}  &= \frac{1}{V}( A_1 \, \tilde H_{1,s} - A_{2} \, \tilde H_{2,s} - A_{3} \, \tilde H_{3,s} ) + \hat J_{sg} - \hat Q_{amb}, \\
    \frac{d \hat u_g}{dt} &= \frac{1}{V}( A_1 \, \tilde H_{1,g} - A_{2} \, \tilde H_{2,g} - A_{3} \, \tilde H_{3,g} ) - \hat J_{sg} - \hat Q_{amb}.
\end{align}
\end{subequations}
The solid-to-gas heat transfer rate $\hat J_{sg}$ is \eqref{eq:solidtogas}. 

\subsection{Transport model}
This section defines the transportation of mass principles in the model. The material fluxes in and out of the cyclone are defined. These directly depend on the gas and solid velocities, which are modeled. The inlet (1) and outlet (2) velocities are pressure drop dependent, while the separation velocity is given explicitly. The amount of solid that is effectively separated depends on the cyclone separation efficiency, which is also modeled in this section.

\subsubsection{Fluxes}
The molar fluxes of the solid and gas streams at the inlet (1), outlet (2), and separation (3) are
\begin{align} 
N_{1, s}&=v_{1} c_{s,in},      & N_{1, g}&=v_{1} c_{g,in},  \\
N_{2, s}&=v_2(1-\eta) c_s,      & N_{2, g}&=v_2 c_g, \\
N_{3, s}&=v_3 \, \eta \, c_s,   & N_{3,g} &= 0.
\end{align}
$\eta$ is the cyclone separation efficiency and $c_{in}$ are the inlet concentrations.

\subsubsection{Inlet and outlet velocities} 
\begin{figure}[tb]
    \centering
    \includegraphics[width=0.3\textwidth]{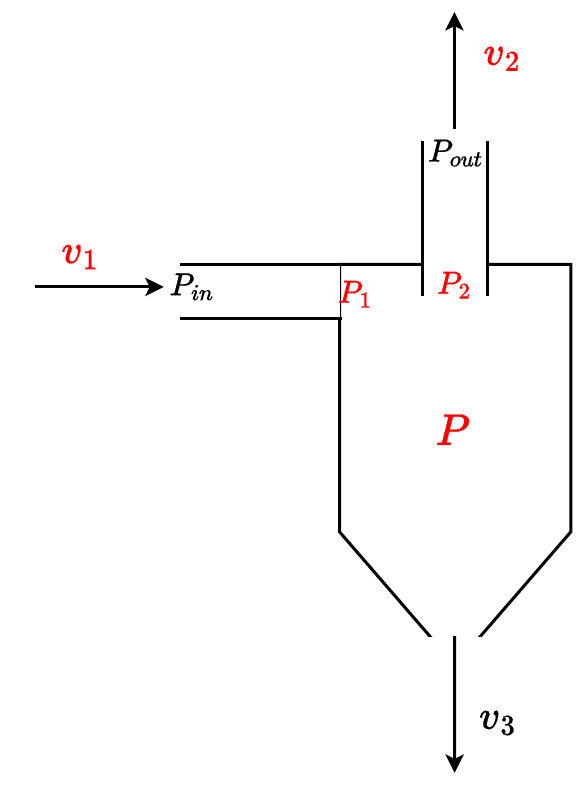}
    \caption{Cyclone pressure points and velocities. The red variables are algebraic variables of the single cyclone model.}
    \label{fig:cyclone_pressure}
\end{figure}

We model the inlet (1) and outlet (2) velocities as function of the pressure drop.
\cite{Chen:2007} present a detailed model that describes the pressure drop in a gas cyclone with tangential inlet. The overall pressure drop in a cyclone is defined as the difference between the pressure at the inlet (1) and the gas outlet (2), i.e.
\begin{equation}
    \Delta P = P_{in} - P_{out}.
\end{equation}
We assume that $P_{in}$ and $P_{out}$ are known (or they are 'external' variables for the cyclone).
Several factors contribute to the pressure drop. These are the expansion loss at the cyclone inlet, $\Delta P_a$, the swirling loss, $\Delta P_b$, the contraction loss at the entrance of the gas outlet (negligible), and the dissipation loss of the gas dynamic energy in the outlet tube, $\Delta P_c$.
\begin{equation}
    \Delta P =\Delta P_a+\Delta P_b+\Delta P_c.
\end{equation}
We define the pressure points $P_1$ and $P_2$ as the pressures after $\Delta P_a$ and $\Delta P_b$, respectively, i.e.
\begin{subequations}
    \begin{align} 
     \Delta P_a & = P_{in} - P_1,\\ 
    \Delta P_b & = P_1 - P_2,\\ 
    \Delta P_c & = P_2 - P_{out}.
    \end{align}
\end{subequations}
Furthermore, we assume that the total pressure in the cyclone, $P$, is an average between $P_1$ and $P_2$, that is
\begin{equation}
    P  = \frac{P_1 + P_2}{2}.
\end{equation}
This is a necessary assumption to relate the total pressure of the cell model to the pressure drops. In this sense, $P$ represents a sort of "average" pressure in the cyclone.
Figure \ref{fig:cyclone_pressure} shows the pressure points and the velocities.

\cite{Chen:2007} derive the following expressions for the first 2 components.
\begin{subequations}\label{eq:dP12}
    \begin{align} 
    \Delta P_a (v_1) & =\left(1-\frac{k_i w_{i n}}{r_c-d_e / 2}\right)^2 \frac{\rho_g v_1^2}{2}, \\
    \Delta P_b (v_1) & =\frac{4 K_A f_0 A_c \tilde{v}_{\theta w}^3}{0.9 \pi D^2 \left( \frac{d_e}{D} \right) ^{1.5 n}}\frac{\rho_g v_1^2}{2}.    \end{align}
\end{subequations}
$\rho_g$ is the density of the gas. This is computed as
\begin{equation}
    \rho_g = \frac{1}{\hat v_g}(M_g \cdot c_g) 
\end{equation}
Notice that dividing by $\hat v_g$ is necessary because of how we defined the concentrations (see \eqref{eq:conc_def}).
The parameters of the expressions are defined as
\begin{subequations}
    \begin{align} 
    k_i & \approx 0.3, \quad K_A=\frac{\pi r_c^2}{h_{i n} w_{i n}}, \quad f_0 = 0.005\left(1+3 \sqrt{c_0}\right), \\ 
    n & =1-\exp \left(-0.26 \,Re^{0.12}\left(1+\left|\frac{h_{e}-h_{i n}}{w_{i n}}\right|\right)^{-0.5}\right), \\ 
    Re & =\frac{\rho_g v_1 2 r_c}{\mu_g K_A \frac{d_e}{D}}, \quad \tilde{r}_c=0.38 \frac{d_e}{D}+0.5 \frac{d_e^2}{D^2}, \\ 
    \tilde{v}_{\theta w} & =\frac{1.11 K_A^{-0.21}\left(\frac{d_e}{D}\right)^{0.16} R e^{0.06}}{\left(1+0.35 {c_0}^{0.27}\right)\left(1+f_0 \frac{A_c}{\pi r_c^2} \sqrt{K_A \frac{d_e}{D}}\right)}.
    \end{align}
\end{subequations}
We refer to \citep{Chen:2007} for further details.
The solid load fraction is 
\begin{equation}
    c_0 = \frac{\sum_{i \in s} M_i \, c_{in,i}}{\sum_{i} M_i \, c_{in,i}},
\end{equation}
where $c_{in}$ are the inlet concentrations.

\cite{Chen:2007} derive a simplified expression for $\Delta P_c$, by using the flow conservation equation, which assumes that the gas volumetric flow rate of the inlet is equal to the one in the outlet. As our model is dynamic, we want to rely on steady state relationships as little as possible. In the outlet tube, most of the gas flows in an annular region close to the wall, while in the vicinity of the center (core region) the axial velocity is close to zero. The pressure loss due to dynamic dissipation depends on the axial velocity in the the annular region, $v_a$, and the mean tangential velocity, $v_\theta$, i.e.
\begin{equation}\label{eq:dP3}
    \Delta P_c = \frac{\rho_g}{2} ( v_\theta ^2 + v_a^2 ).
\end{equation}
$v_a$ and $v_\theta$ can be expressed as a function of the average outlet velocity $v_2$ \citep{Hoffmann:2002}
\begin{equation}\label{eq:a_theta}
    v_a = \frac{v_2}{1 - R_{cx}^2}, \quad v_\theta = \sqrt{\frac{2 R_{cx}^3}{(1 - R_{cx}^2)^3}} v_2,
\end{equation}
where
\begin{equation}
    R_{cx}  = \frac{\tilde r_c  D}{d_e}.
\end{equation}
Inserting \eqref{eq:a_theta} in \eqref{eq:dP3} gets
\begin{equation}
    \Delta P_c (v_2) = \left(\frac{2 R_{cx}^3 - R_{cx}^2 - 1}{(1- R_{cx}^2)^3}\right) \frac{\rho_g v_2^2}{2}.
\end{equation}

We now have expressions of $\Delta P_a$ and $\Delta P_b$ as a function of the inlet velocity $v_1$ \eqref{eq:dP12}, and an expression of $\Delta P_c$ as a function of $v_2$ \eqref{eq:dP3}.

\subsection{Separation velocity}
\cite{Hoffmann:2002} compute the separation velocity, $v_3$, or 'wall axial velocity', as
\begin{equation}
    v_3 = \frac{0.9 A_1 v_1}{\pi  (r_c^2 - r_m^2)}, \quad r_m = \sqrt{r_c \frac{d_e}{2}}.
\end{equation}

\subsubsection{Separation efficiency}
Not all the suspended solid particles in the gas get separated by the cyclone. The amount of solid particles that is successfully separated from the gas phase depends on the separation efficiency of the cyclone. \cite{Muschelknautz:1993} formulate the following model for the separation efficiency
\begin{equation}
    \eta=1-\frac{c_{0 L}}{c_0}+\frac{c_{0 L}}{c_0} \exp \left(-\left(\frac{d^*}{d_A}\right)^n\right),
\end{equation}
where
\begin{subequations}
\begin{align}    
n & = 1.25, \\
c_{0L} & =0.025\left(\frac{d^*}{d_{\text {med }}}\right)\left(10 c_0\right)^k, \\
d_A & =\frac{d_e^*}{0.7^{\frac{1}{n}}}, \\
k & = \begin{cases}
    -0.11-0.10 \ln (c_0) & c_0 \geq 0.1,\\
    0.15 & c_0 < 0.1.
\end{cases}
\end{align}
\end{subequations}
$c_{0L}$ is the solid loading limit. The particle cut-size, $d^*$, and suspended characteristic particle size, $d_e^*$, are given by
\begin{equation}
     d^*  =\sqrt{\frac{18 \cdot 0.9 \, \mu \, A_{1} v_{1}}{2 \pi (h_t - h_e) u(r_i)^2 \Delta \rho }}, \quad d_e^*  =\sqrt{\frac{18 \mu}{\Delta \rho \, \bar{z}} w_{s 50}}.
\end{equation}
The parameters are
\begin{subequations}
\begin{align}
& w_{s 50} =\frac{0.5 \cdot 0.9 A_{1} v_{1}}{A_w}, \quad \Delta \rho = \rho_s-\rho_g,  \\
&A_w =2 \pi r_c\left(h_t-h_c\right)+\pi\left(r_c+r_2\right) \sqrt{\left(r_c-r_2\right)^2+\frac{h_c^2}{4}}, \\
&r_i=\frac{d_e}{2}, \quad u(r)=\frac{u_c \frac{r_c}{r}}{1+f_0 \frac{A_w u_c}{2 A_{1} v_{1}} \sqrt{\frac{r_c}{r}}}, \\
&\bar{z} =\frac{u(r_1) u(r_2)}{\sqrt{r_2 r_{1}}}, \quad  r_1 = r_c - \frac{w_{in}}{2}, \\
&r_2 =r_c-\frac{1}{2}\left(r_c-\frac{d_d}{2}\right), \quad u_c =\frac{r_{1}}{\alpha r_c} v_{1}, \quad \beta=\frac{w_{i n}}{r_c} \\
& \alpha =\frac{1}{\beta}\left(1-\sqrt{1+(\beta^2 - 2\beta) \sqrt{1-\frac{\beta(2-\beta)\left(1-\beta^2\right)}{1+c_0}}}\right).
    \end{align}
\end{subequations}
The density of the solid is
\begin{equation}
    \rho_s = \frac{1}{\hat v_s}(M_s \cdot c_s) 
\end{equation}
Notice that this formulation of the Muschelknautz model does not use the whole particle size distribution, but only the median size.

\subsubsection{Algebraic relations}
Some algebraic relations are needed to complete the model. These are
\begin{subequations}
\label{eq:AlgCyclone}
\begin{align}
    & V(T_{s},P,c_{s}) + V(T_{g},P,c_{g}) - 1 = 0, \label{eq:AlgCyclone1} \\ 
    & U(T_{s},P,c_{s}) - \hat u_{s} =0, \label{eq:AlgCyclone2}\\
    & U(T_{g},P,c_{g}) - \hat u_{g} =0, \label{eq:AlgCyclone3}\\
    & \Delta P_a(v_1) - (P_{in} - P_1) = 0, \label{eq:AlgCyclone4}\\ 
    &\Delta P_b(v_1) - (P_1 - P_2) = 0, \label{eq:AlgCyclone5}\\ 
    &\Delta P_c(v_2) - (P_2 - P_{out}) = 0, \label{eq:AlgCyclone6}\\
    &P - (P_1 + P_2)/2 = 0. \label{eq:AlgCyclone7}
    \end{align}
\end{subequations}
\eqref{eq:AlgCyclone1} requires that the volume of the solid and the gas in the cyclone adds up to the volume of the cyclone. \eqref{eq:AlgCyclone2} and \eqref{eq:AlgCyclone3} require that the differential variables $\hat u_{s}$ and $\hat u_{g}$ match the quantities computed via the thermodynamic functions. \labelcref{eq:AlgCyclone4,eq:AlgCyclone5,eq:AlgCyclone6,eq:AlgCyclone7} introduce the pressure drop model.
\subsubsection{Summary}
Each cyclone is modeled by a system of DAEs, consisting of the differential equations \eqref{eq:cyclonemassbalance} and \eqref{eq:cycloneenergybalance}, and the algebraic equations \eqref{eq:AlgCyclone}. For each cyclone, the differential and algebraic variables are
\begin{equation}
    x^{Cyc} = \begin{bmatrix}
        c \\ \hat u_s \\ \hat u_g
    \end{bmatrix}, \quad y^{Cyc} = \begin{bmatrix}
        T_s \\ T_g \\ P \\ P_1 \\ P_2 \\ v_1 \\ v_2
    \end{bmatrix}
\end{equation}
$c$ is intended as a vector like usual.

\section{Other components}

\subsection{Filter}
For simplicity, we assume that the dust is completely removed from the stream. This means that we assume that the filter has 100\% efficiency (ideal filter).

\subsection{Circulation fan}
The fan may be modeled by the static relation
\begin{equation}\label{eq:fan}
    \eta_{fan} = \frac{F_{vol} \Delta P}{P_{fan}}.
\end{equation}
$\eta_{fan}$ is the efficiency of the fan.
$F_{vol}$ is the volumetric flow rate through the fan, $P_{fan}$ is the power used by the fan. The pressure increase is
\begin{equation}
    \Delta P = P_{out}^{fan} - P_{in}^{fan}.
\end{equation}

\subsection{Gas purge}
The gas purge before the recirculation is modeled as a stream splitter, whose static mass balance reads
\begin{equation}
    f_{in} = \alpha_{purge} f_{in} + f_{recirc}.
\end{equation}
$f_{in}$ is the flow rate of the incoming stream, $\alpha_{purge} \in [0,1]$ is the percentage of purge, and $f_{recirc}$ is the recirculated flow. 

\subsection{Fresh air mixer}
The fresh air is mixed to the remaining recirculated gas. This is modeled by a stream mixer, whose static mass and energy balances read 
\begin{subequations}
    \begin{align}
        f_{mix} & = f_{recirc} + f_{fresh}, \\
        \begin{split}
            H(T_{mix},P_{out}^{fan},f_{mix}) & = H(T_{g}^{cyc1},P_{out}^{fan},f_{recirc}) \\
            & + H(T_{fresh},P_{out}^{fan},f_{fresh}). \label{eq:enbalance_mixer}
        \end{split}
    \end{align}
\end{subequations}
The energy balance is necessary because the streams mix with different temperatures.
$f_{fresh}$ is the fresh air feed flow rate, at temperature $T_{fresh}$, and $f_{mix}$ is the mixed stream flow rate, at temperature $T_{mix}$. We assume that the pressure is kept constant. Notice that the recirculated gas has temperature $T_{g}^{cyc1}$, i.e. the gas temperature when leaving cyclone 1. This is due to the assumption of ideal, adiabatic connecting tubes.

\subsection{Electric hot gas generator}
The electric hot gas generator is a device that deploys resistive heating as mean of energy transfer between electrical current and the working fluid. The electrical energy is converted to heat via the Joule effect. The static energy balance of the gas stream that undergoes heating is
\begin{equation}\label{eq:EHGG}
    H(T_{mix},P_{out}^{fan},f_{mix}) + P_{EHGG} = H(T_{g,in}^{Calc},P_{out}^{fan},f_{mix}),
\end{equation}
where $P_{EHGG}$ is the electrical power.

\section{Plant-wide model and connection of the units}\label{sec:Connection}
\begin{figure}[tb]
    \centering
    \includegraphics[width=0.5\textwidth]{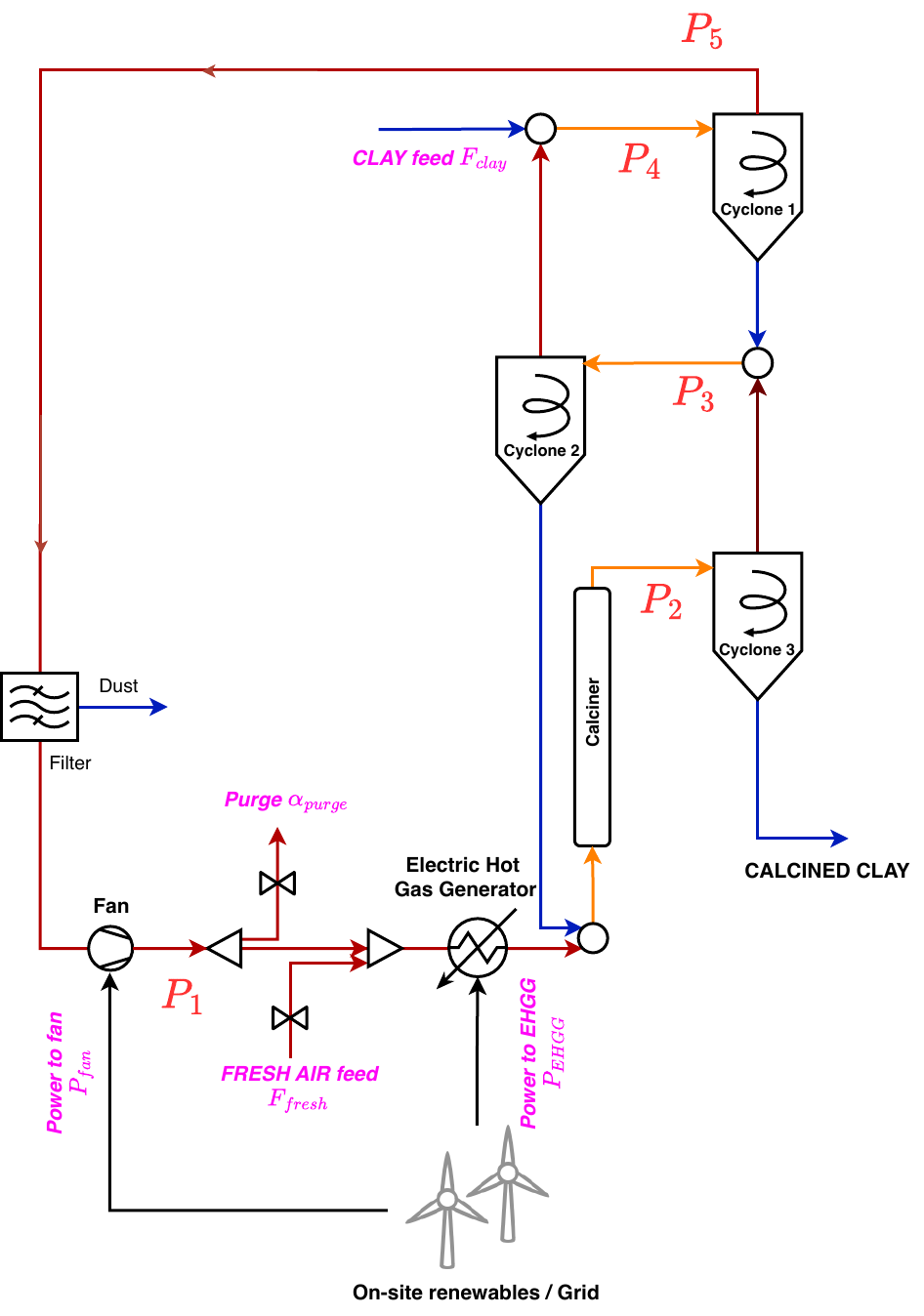}
    \caption{Pressure points in the plants.}
    \label{fig:pressurepoints}
\end{figure}

The connection between the units can be modeled in several ways, depending on the degree of accuracy of the model that is required. The most precise way would be to model every tube as a plug-flow reactor (as for the calciner), and therefore consider reaction, transport and heat transfer phenomena also in the connections. This would result in a large increase of the number of system states. To avoid over-complicating the model, we propose a simpler approach.
We assume that the pressure at the outlet of one unit equals the pressure at the inlet of the next unit. This assumption is good enough when the connecting tubes are short, and therefore the pressure drop inside the tubes is negligible. By using this approach, five pressure points are identified (see Figure \ref{fig:pressurepoints}). For each unit, we have that
\begin{subequations}
    \begin{align}
    & P_{in}^{Calc} = P_1, \quad P_{out}^{Calc} = P_2 \\
    & P_{in}^{Cyc3} = P_2, \quad P_{out}^{Cyc3} = P_3 \\
    & P_{in}^{Cyc2} = P_3, \quad P_{out}^{Cyc2} = P_4 \\
    & P_{in}^{Cyc1} = P_4, \quad P_{out}^{Cyc3} = P_5 
\end{align}
\end{subequations}
The pressure points are algebraic variables in the model. Five extra equations are necessary to solve the model. The flow continuity equations between cyclones and calciner may be used to this scope
\begin{subequations}\label{eq:continuity}
\begin{align}
    &A_1  v_1^{Cyc2} = A_2 v_2^{Cyc3}, \\
    &A_1  v_1^{Cyc1} = A_2 v_2^{Cyc2}, \\
    &A_{Calc} v_{N_z+1/2}^{Calc} = A_1 v_1^{Cyc3}.
\end{align}
\end{subequations}
Moreover, the volumetric flow rate that enters the calciner must satisfy
\begin{equation}\label{eq:calc_alg}
    A_{Calc} v_{1/2}^{Calc} = V_g(T_{g,in}^{Calc},P_1,f_{mix}) + V_s(T_{s}^{Cyc2},P_1,f_{3,s}^{Cyc2})
\end{equation}
The last algebraic equation is given by the fan equation \eqref{eq:fan}, where the pressure increase is $\Delta P = P_1 - P_5$. Note that \eqref{eq:continuity} do not need to be necessarily  implemented as algebraic equations, since one of the velocities in each equations can be explicitly isolated and therefore removed as algebraic variable.

The transport of the material from one unit to another is realized by appropriately selecting the inlet concentrations of each unit. These can be calculated by using the assumption that the molar flow rate at the outlet of one unit should be equal to the one at the inlet of the next unit. We remark that a molar flow rate can be expressed as $f = A \, v \, c$. The inlet temperatures of solid and gas are also easily carried from the previous unit.
Notice that cyclone 1 is different from the other two, since it receives a feed of fresh clay from the outside. Its solid inlet concentration is
\begin{equation}
    c_{s,in}^{Cyc1} = \frac{f_{clay}}{v_1^{Cyc1} A_1} + (1- \eta^{Cyc2}) c_s^{Cyc2}.
\end{equation}
$f_{clay}$ is the molar flow rate of fresh clay entering the cyclone, while the second term is the unseparated solid in the gas coming from cyclone 2. The temperature of the fresh clay is $T_{clay}$.

\begin{rem}  Another solution to the connection problem is to neglect transport delay, heat transfer and reaction, but to consider the pressure drop. The connection between the units may be realised by lumping the flow resistance in the connecting tube as a pressure node \citep{Hansen:1998}. The flow rate $F$ between 2 units (1 and 2) may be expressed as
\begin{equation}\nonumber
    F = C \sqrt{\frac{(P_1 - P_2)(P_1+P_2)}{T_1}},
\end{equation}
where $C$ represents the flow resistance.
\end{rem}
\subsection{Summary of the full model}
The full model is obtained by combining the equations of all the units and connecting them as explained in the previous subsection. The equations consist of \eqref{eq:massODE}, \eqref{eq:energyODE}, \eqref{eq:algebraic} for the calciner, \eqref{eq:cyclonemassbalance}, \eqref{eq:cycloneenergybalance}, \eqref{eq:AlgCyclone} for each cyclone, the secondary units \eqref{eq:fan}, \eqref{eq:enbalance_mixer}, \eqref{eq:EHGG}, and the connecting equations \eqref{eq:continuity}, \eqref{eq:calc_alg}. The differential and algebraic variables of the full model are
\begin{equation}
    x = \begin{bmatrix}
        x^{Cyc1} \\ x^{Cyc2} \\ x^{Cyc3} \\ x^{Calc}
    \end{bmatrix}, \quad 
    y = \begin{bmatrix}
        y^{Cyc1} \\ y^{Cyc2} \\ y^{Cyc3} \\ y^{Calc} \\ P_1 \\ P_2 \\ P_3 \\ P_4 \\ P_5 \\ T_{mix} \\ T_{g,in}^{Calc}
    \end{bmatrix}.
\end{equation}
The manipulated variables are
\begin{equation}
    u = \begin{bmatrix}
        f_{clay} \\
        \alpha_{purge} \\
        f_{fresh} \\
        P_{fan} \\
        P_{EHGG}
    \end{bmatrix}.
\end{equation}
Notice that, depending on the control structure of the whole cement plant, $P_{EHGG}$ might be a disturbance variable for the clay plant rather than being directly controllable.
The disturbance variables are
\begin{equation}
    d = \begin{bmatrix}
        d_{med}\\
        T_{fresh} \\
        T_{clay}
    \end{bmatrix}.
\end{equation}
These can be considered parameters rather than disturbances if assumed constant.
The model outputs of interest are the calcined clay production rate and its degree of calcination. These can be computed as
\begin{equation}
    z = \begin{bmatrix}
        CC \\ CD
    \end{bmatrix} = \begin{bmatrix}
        A_3 v_3^{Cyc3} \eta^{Cyc3} (c_{s}^{Cyc3} \cdot M_s) \\
         \frac{c_{A}^{Cyc3} M_A} {c_{AB_2}^{Cyc3} M_{AB_2} + c_{A}^{Cyc3} M_A}
    \end{bmatrix}.
\end{equation}

\subsection{Model derivatives}
\begin{figure*}
    \centering
    \includegraphics[width=0.9\linewidth]{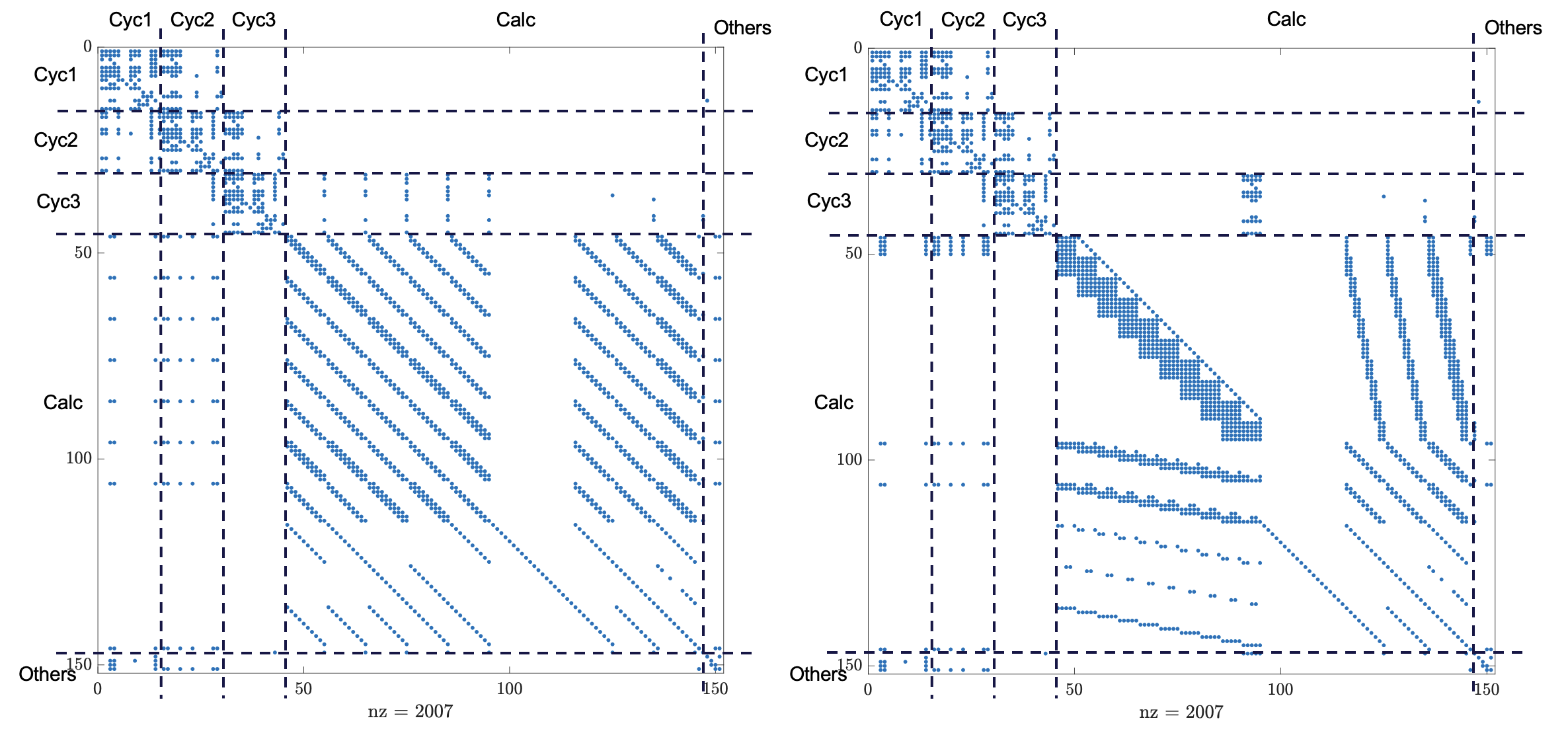}
    \caption{Non-zero elements (spy) of the Jacobian of the full model, $\frac{\partial F}{\partial z}$. Two different variable orderings, by chemical species, on the left, and by cell, on the right, are considered for the calciner.}
    \label{fig:Jacobian}
\end{figure*}

In order to simulate the model, the Jacobian of the model, 
\begin{equation}
    J = \frac{\partial F}{\partial z},
\end{equation}
is required. $F$ are the combined differential and algebraic equations, and $z$ are the combined differential and algebraic variables (notice that these are not necessarily grouped). By using a DAE solver like \verb|ode15s| in Matlab, the derivatives can be approximated by finite difference. Nonetheless, computing the Jacobian analytically, and providing it to the solver, can improve significantly the performance. This is especially relevant, if not crucial, when performing optimization using the model (which is, however, outside of the scope of this paper). We compute the derivatives analytically, and verify them by comparing them against the finite difference approximation. Figure \ref{fig:Jacobian} shows the non-zero elements (spy) of the Jacobian matrix, $J$. We consider and show two different variable ordering for the calciner, i.e. by chemical species and by cell. We find that using the cell ordering leads to slightly better performance of the solver. Overall, we observe a 6x speed-up when providing the analytical Jacobian to \verb|ode15s|, compared to the default finite difference approximation.
\section{Simulation results}\label{sec:simulation}
In this section, we demonstrate the use and the implementation of the dynamic model by executing a simulation with changing input. The model is implemented in Matlab and simulated using \verb|ode15s|. The algebraic equations are specified by defining a mass matrix. The units are connected by using the 5 pressure points, as described in Section \ref{sec:Connection}.

\begin{figure}[tb]
    \centering
    \includegraphics[width=0.45\textwidth]{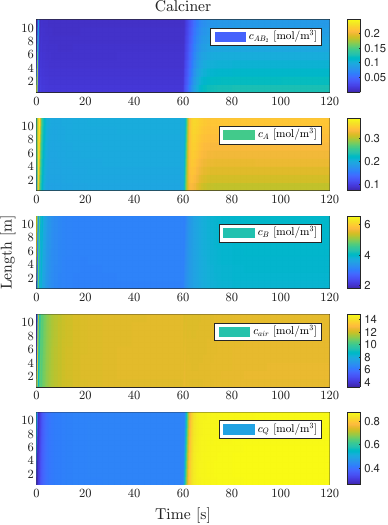}
    \caption{Concentrations in space and time in the calciner.}
    \label{fig:calciner_conc}
\end{figure}

\begin{figure}[tb]
    \centering
    \includegraphics[width=0.45\textwidth]{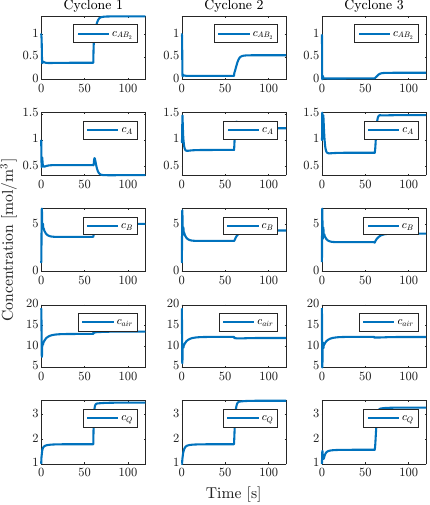}
    \caption{Concentrations in the cyclones.}
    \label{fig:cyclones_conc}
\end{figure}

\begin{figure}[tb]
    \centering
    \includegraphics[width=0.45\textwidth]{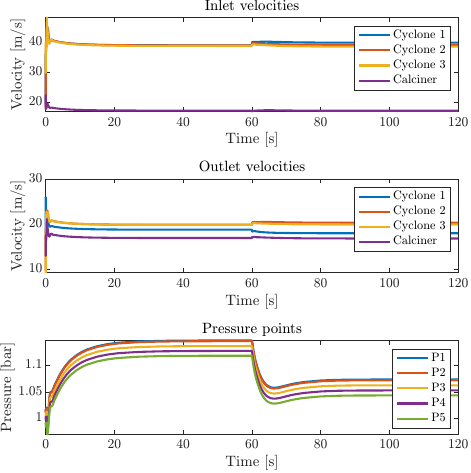}
    \caption{Inlet and outlet velocities of the calciner and the cyclones, and the pressure points.}
    \label{fig:velocities}
\end{figure}

\begin{figure}[tb]
    \centering
    \includegraphics[width=0.45\textwidth]{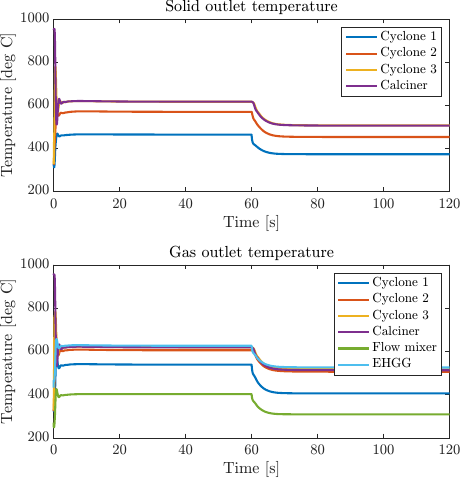}
    \caption{Temperature of the solid and the gas at the outlet of each unit.}
    \label{fig:temperatures}
\end{figure}

\begin{figure}[tb]
    \centering
    \includegraphics[width=0.45\textwidth]{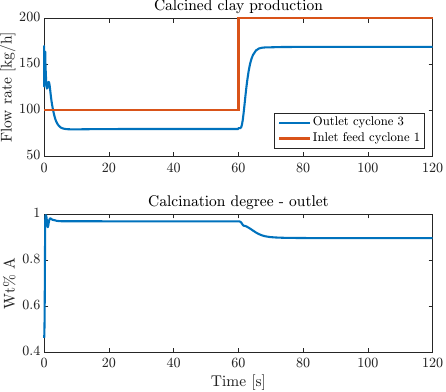}
    \caption{Outputs of the system: production of calcined clay and its calcination degree. The inlet clay feed is also plotted.}
    \label{fig:Outputs_sim}
\end{figure}

\begin{figure}[tb]
    \centering
    \includegraphics[width=0.5\textwidth]{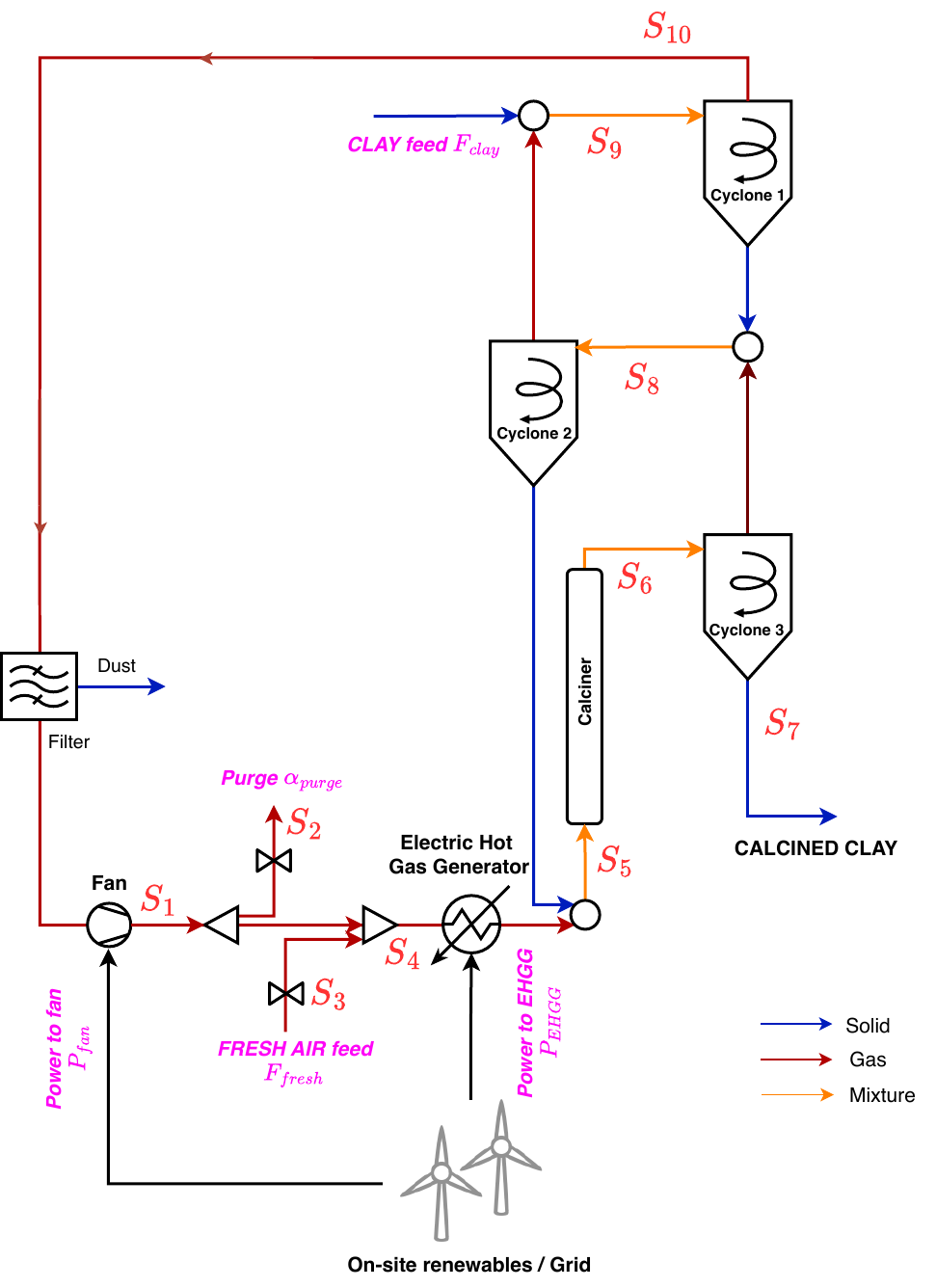}
    \caption{Relevant material streams numbering.}
    \label{fig:Streams}
\end{figure}

For the cyclones, the Stairmand design (with slot inlet) is chosen. Its geometric parameters are reported in Table \ref{tab:Stairmand}. For the calciner, the length is set to $L = 12$ m, the diameter to $d = 0.18$ m, and the diffusion parameters to $D=0.1$. We assume that all reactors are adiabatic, i.e. $Q_{amb}^i = 0$ for $i=\{Calc,Cyc1,Cyc2,Cyc3\}$. The median particle size and the solid-to-gas heat transfer coefficient are
\begin{equation}
    d_{med} = 7.61 \times 10^{-6} \text{m}, \qquad k_{sg} = 200 \mathrm{\frac{J}{s\, m^2 K}}.
\end{equation}
The fresh clay is assumed to be 68\% kaolinite and 32\% quartz in weight. The efficiency of the circulation fan is assumed to be $\eta_{fan}=0.8$. We simulate the system for 2 min. We conduct a step experiment by doubling the clay feed, to investigate the effect on the system and the output. We apply a clay feed of 100 kg/h for 1 min and 200 kg/h for the other minute. The other inputs are kept constant to 
\begin{equation}
    P_{fan} = 1.2 \text{ kW}, \quad P_{EHGG} = 40  \text{ kW}, \quad \alpha_{purge} = 0.3.
\end{equation}
The fresh air feed is set to 150 kg/h of air, of which 10\% is water vapor. The results are now presented. Figure \ref{fig:calciner_conc} shows the concentrations in space and time in the calciner. Figure \ref{fig:cyclones_conc} shows the concentrations in the three cyclones. Figure \ref{fig:velocities} plots the inlet and outlet velocities of the calciner and the cyclones, and the pressure points. Figure \ref{fig:temperatures} shows the outlet temperature of each unit, of both the solid and the gas phase. Figure \ref{fig:Outputs_sim} finally shows the outputs of the system, i.e. the production of calcined clay and its calcination degree. The step input of clay feed is also plotted. The effect of doubling the clay feed on the system is significant. The overall pressure in the system drops, as well as the temperatures. This is of course due to the fact that the solid now needs more energy to be heated up. In the first part of the simulation, the calcination degree reaches $\sim$97\% at steady state, while in the second part it drops to $\sim$90\%. The calcination drops because more energy is necessary. The system reaches a new steady state within 20 sec. Figure \ref{fig:Streams} defines the material streams $S = \{S_1,...,S_{10}\}$. Table \ref{tab:streams} reports relevant numerical values at the last steady state of the simulation ($t=120$ s). We report the solid mass flow rate, $F_s$, the gas mass flow rate, $F_g$, the pressure, $P$, the solid temperature, $T_s$, the gas temperature, $T_g$. Moreover, we report the metakaolin weight percentage content (or calcination degree), $w_A$, and the water weight percentage content, $w_B$. These are defined as
\begin{equation}
    w_A = \frac{F_A}{F_{AB_2}+F_A}, \quad w_B = \frac{F_B}{F_{air}+F_B}= \frac{F_B}{F_g}.
\end{equation}
By looking at these statistics, we can follow the calcination process through the different units. After the pre-heating cyclones, about 66\% conversion is achieved. The calcination is then completed in the calciner ($\sim 78$\%) and the separation cyclone ($\sim 90\%$).  The experiment confirms that the system dynamics are non-trivial, and the model can enable process optimization and model-based control.
\begin{table}[tb]
\centering
\label{tab:Stairmand}
\caption{Geometric parameters of the Stairmand cyclone (cyclone diameter of $D=0.3$ m)}
\begin{tabular}{|l|l|}
\hline
Parameter & Value\\
\hline
$h_{in}/D$ & 0.5 \\
$w_{in}/D$ & 0.2 \\
$d_e/D$ & 0.5 \\
$h_e/D$ & 0.5 \\
$h_t/D$ & 4 \\
$h_c/D$ & 2.5 \\
$d_d/D$ & 0.37 \\
\hline
\end{tabular}
\end{table}

\begin{table*}[tb]
    \centering
    \caption{Selected material streams at steady state ($t=120$ s). For each stream, the solid and gas mass flow rates, pressure, temperatures and percentage of metakaolin and water are reported.}
    \begin{tabular}{|c|c|c|c|c|c|c|c|}
\hline
Stream & $F_s$ [kg/h] & $F_g$ [kg/h] & $P$ [bar] & $T_s$ [${}^\circ C$] & $T_g$ [${}^\circ C$] & $w_A$ [\%] & $w_B$ [\%] \\
 \hline
$S_1$ & 0 & 554.14 & 1.0739 & - & 405.35 & - & 18.793 \\
$S_2$ & 0 & 166.24 & 1.0739 & - & 405.35 & - & 18.793 \\
$S_3$ & 0 & 150 & 1.0739 & - & 30 & - & 10 \\
$S_4$ & 0 & 537.9 & 1.0739 & - & 308.81 & - & 16.341 \\
$S_5$ & 190.75 & 537.9 & 1.0739 & 452.86 & 525.55 & 66.287 & 16.341 \\
$S_6$ & 188.55 & 540.1 & 1.0724 & 505.75 & 512.82 & 78.476 & 16.682 \\
$S_7$ & 169.12 & 0 & - & 506.56 & - & 89.635 & - \\
$S_8$ & 218.04 & 542.04 & 1.0629 & 372.77 & 516.99 & 22.27 & 16.981 \\
$S_9$ & 217.85 & 551.49 & 1.0532 & 68 & 505.77 & 5.2677 & 18.403 \\
$S_{10}$ & 14.645 & 554.14 & 1.0437 & 372.77 & 405.35 & 16.61 & 18.793 \\
\hline
\end{tabular}
    \label{tab:streams}
\end{table*}

\section{Conclusion}\label{sec:conclusion}
A dynamic model of an electric clay calcination plant is presented, discussed and simulated. The model consists of several building blocks, which are formulated in a way that allows easy modifications, if necessary. For example, extra side reactions can be added, and different transport models can be investigated, with minimum effort in the implementation. The rigorous thermodynamic functions allow realistic evaluation of heat transfer phenomena at changing conditions. The model describes in details the dynamics in the main units of the plant, i.e. the cyclones and the calciner. These units are then connected in a loop, and a system of DAEs arises. The plant-wide model can therefore be easily simulated with a standard DAE solver. Finally, we perform and present a simulation where we apply a step input in the clay feed and observe the dynamics in system.

The dynamic model is an efficient proxy to perform easy and quick simulations with different process inputs and conditions. The model can unlock process control and optimization of the plant. By linearizing the model, linear MPC can be developed. Moreover, PID controllers can be tuned and their performance can be evaluated by using the model. Furthermore, steady state and cost optimization can be also performed. The production cost can be optimized, by minimizing the use of energy, while maximizing the production rate.
\section*{Acknowledgement}
This project has been partly funded by the Danish Energy Technology and Demonstration Program (EUDP) under the Danish Energy Agency in the EcoClay project no. 64021-7009.



\bibliographystyle{elsarticle-harv} 


\bibliography{cas-refs}

\end{document}